\documentstyle[11pt]{article}
\begin{document}

\title{{\bf Differential Forms and the Noncommutative Residue for Manifolds with Boundary\\
 in the Non-product Case}
\thanks{partially supported by CNSF No. 10231010 and NSF of Zhejiang
 Province, No. 101037 and Science Foundation for Young Teachers of Northeast Normal University (No. 20060102)
 }}
\author{ Yong Wang\\
\thanks{also at School of Mathematics and Statistics, Northeast Normal University}
{\scriptsize \it  Center of Mathematical Sciences Zhejiang
University Hangzhou Zhejiang 310027, China ;}\\
{\scriptsize \it E-mail: wangy581@nenu.edu.cn}}

\date{}
\maketitle

\noindent {\bf Abstract}~~ In this paper, for an even dimensional
compact manifold with boundary which has the non-product metric near
the boundary, we use the noncommutative residue to define a
conformal invariant pair. For a $4$-dimensional manifold, we compute
this conformal invariant pair under some conditions and point out
the way of computations in the general.
\\
\noindent{\bf Subj. Class.:}\quad Noncommutative differential
geometry; Noncommutative global analysis.
 Noncommutative residue for manifolds with
boundary; non-product metric; conformal invariant. \\

\section{Introduction}
\quad Since the noncommutative residue was found in
[Ad],[M],[Gu],[Wo], it was applied to many branches of mathematics.
Especially, it was as the noncommutative counterpart of the integral
in NCG by [C1]. The noncommutative residue also had been used to
derive the gravitational action in the framework of NCG in [K],
[KW]. In [C2], Connes used the noncommutative residue to find a
conformal 4-dimensional Polyakov action analogy. In [U],
Connes' result was generalized to the higher dimensional case.\\
\indent The noncommutative residue on Boutet de Monvel algebra for
manifolds with boundary was found in [FGLS]. In [S], Schrohe gave
the relation between the Dixmier trace and the noncommutative
residue for manifolds with boundary. In [Wa1], the author proved a
Kastler-Kalau-Walze type theorem for manifolds with boundary and for
the boundary flat case, he gave two kinds of operator theoretic
explaination of the gravitational action on boundary. In [Wa2], the
author generalized the results in [C2] and [U] to the case of
manifolds with boundary which have a product metric near the
boundary. A natural question is to define and compute a conformal
invariant pair in the non-product metric case. In this paper, for an
even dimensional compact manifold with boundary which has a
non-product metric near the boundary, we define a conformal
invariant pair. When $n=4$, we compute this conformal invariant pair
under some conditions and point out the way of computations in the
general. As a corollary, when $n=4$, for some special non-product
metrics, we get the conformal invariant on the
boundary vanishes which generalizes partially a result in [Wa2]. \\
\indent This paper is organized as follows: In Section 2, we define
a conformal invariant pair associated to an even dimensional compact
manifold with boundary which has a non-product metric near the
boundary. In Section 3, for a $4$-dimensional manifold, we compute
this conformal
invariant pair under some conditions. Some remarks on computations in the general case when $n=4$ will be given in Section 4.\\

\section{The Conformal Invariant Pair
$(\Omega_n(f_1,f_2),\Omega_{n-1}(f_1,f_2))$}

\quad Let $M$ be an even dimensional compact oriented Riemaniann
manifold with boundary $\partial M$ and $U\subset M$ be the collar
neighborhood of $\partial M$ which is diffeomorphic to $\partial
M\times [0,1)$. Write ${\rm dim}M=n$. Let $g^M$ be the metric on $M$
which has the following form on $U$
$$ g^M=\frac{1}{h(x_n)}g^{\partial M}+dx_n^2,\eqno(2.1)$$
where $g^{\partial M}$ is the metric on  ${\partial M}$; $h(x_n)\in
C^{\infty}([0,1))=\{g|_{[0,1)}|g\in C^{\infty}((-\varepsilon,1))\}$
for some $\varepsilon>0$ and satisfies $h(x_n)>0,~h(0)=1$ where
$x_n$ denotes the normal directional coordinate.\\
    \indent In this section, we will construct a conformal invariant pair
$(\Omega_n(f_1,f_2),\Omega_{n-1}(f_1,f_2))$ associated to $M$. The
fundamental setup is the same as Section 2 and Section 3 in [Wa2].
Recall that in Section 4 of [Wa2], we consider the product metric
case, i.e. $h(x_n)\equiv 1$. We can use a canonical way to construct
a metric $\widetilde{g}$ on the double manifold
$\widehat{M}=M\cup_{\partial M} M$ through taking $\widetilde{g}=g$
on both copies of $M$, then $\widetilde{g}$ is well defined by
$h=1$. But for the general $h$, this is not correct. So we need to
use another way
to construct a conformal invariant pair associated to $M$.\\
\indent   By the definition of $C^{\infty}([0,1))$ and $h>0$,  there
exists $\widetilde{h}\in C^{\infty}((-\varepsilon,1))$ such that
$\widetilde{h}|_{[0,1)}=h$ and $\widetilde{h}>0$ for some
sufficiently small $\varepsilon>0$. Using partition of unity
Theorem, then there exists a metric $\widehat{g}$ on $\widehat{M}$
which has the form on $U\cup_{\partial M}\partial M\times
(-\varepsilon,0]$
$$ g^M=\frac{1}{\widetilde{h}(x_n)}g^{\partial M}+dx_n^2,\eqno(2.2)$$
such that $\widehat{g}|_M=g.$ Nextly we fix a metric $\widehat{g}$
on the $\widehat{M}$ such that $\widehat{g}|_M=g$. Denote by
$[(M,g)]$ a conformal manifold. The way of constructing a conformal
invariant pair associated to $[(M,g)]$ is as follows. As in [C2] or
[U], we consider the following operator on the
manifold $(\widehat{M},\widehat{g})$,\\
$$F_{\widehat{g}}:=\frac{d\delta-\delta d}{d\delta+\delta
d}:\wedge^{\frac{n}{2}}(T^*\widehat{M})\rightarrow\wedge^{\frac{n}{2}}(T^*\widehat{M}),\eqno(2.3)$$
then $F_{\widehat{g}}$ does not depend on the choice of the metric
in the conformal class $[(\widehat{M},\widehat{g})]$.
 Now similar to (3.5) and (3.6)
in [Wa2], for $f_0,f_1,f_2\in C^{\infty}(M)$ and $f_0$ not depending
on $x_n$ near the boundary, we define the form pair
$(\Omega_n(f_1,f_2)(\widehat{g}),\Omega_{n-1}(f_1,f_2)(\widehat{g}))$
through the following equality:\\
$~~\\ \widetilde{{\rm
Wres}}(\pi^+f_0[\pi^{+}F_{\widehat{g}},\pi^+f_1][\pi^+F_{\widehat{g}},\pi^+f_2])$\\
$$~~~~~ =\int _{M}f_0\Omega_n(f_1,f_2)(\widehat{g})+\int _{\partial M}f_0|_{\partial M}
\Omega_{n-1}(f_1,f_2)(\widehat{g}).\eqno(2.4)$$ \noindent By the
definition of $\pi^+F_{\widehat{g}}$ in the Boutet de Monvel
algebra, the left term of (2.4) is well defined. We hope to
generalize the results in [C2] and
[U], so as in [U], we take\\
\noindent $\Omega_{n}(f_1,f_2)(\widehat{g})=$\\
$$\int_{|\xi|=1}{\rm
tr}\left[
\sum\frac{1}{\alpha'!\alpha''!\beta!\delta!}D^{\beta}_x\bar{f_1}
D^{\alpha''+\delta}_x\bar{f_2}
\partial^{\alpha'+\alpha''+\beta}_{\xi}\sigma^{F_{\widehat{g}}}_{-j}
\partial^{\delta}_{\xi}D^{\alpha'}_x\sigma^{F_{\widehat{g}}}_{-k}\right]\sigma(\xi)d^nx\left|_M\right.
,\eqno(2.5)$$
 \noindent where $\sigma_{-j}^{F_{\widehat{g}}}$ denotes the order
$-j$ symbol of ${F_{\widehat{g}}}$; $\overline{f_1},~\overline{f_2}$
are the extensions to $\widehat{M}$ of $f_1,f_2$,
$D^\beta_x=(-i)^{|\beta|}\partial^\beta_x$, and the sum is taken
over $|\alpha'|+|\alpha''|+|\beta|+|\delta|+j+k=n;
|\beta|\geq1,|\delta|\geq1; \alpha',\alpha'',\beta,\delta \in {\bf
Z}^n_+; j,k\in {\bf Z}_+.$ Then $\Omega_{n}(f_1,f_2)(\widehat{g})$
does not depend on the extensions of $f_1,f_2$. By Theorem 3.1 and
(3.19) in [Wa2], then $\Omega_{n-1}(f_1,f_2)(\widehat{g})$ is
uniquely determined by (2.4), (2.5) as follows:
$$\Omega_{n-1}(f_1,f_2)(\widehat{g})=\sum^{\infty}_{j,k=0}\sum
\sum^{-r}_{|\beta|=1}\sum^{-l}_{|\delta|=1}
\frac{(-i)^{j+k+1+|\alpha|+|\beta|+|\delta|}}{\alpha!\beta!\delta!(j+k+1)!}$$
$$\times
\int_{|\xi'|=1} \int^{+\infty}_{-\infty}{\rm trace}\left\{
\partial^j_{x_n}\left[\partial^\beta_xf_1\partial^\alpha_{\xi'}\partial^k_{\xi_n}\pi^+_{\xi_n}\partial^\beta_{\xi}
\sigma_{r+|\beta|}^{F_{\widehat{g}}}\right]|_{x_n=0}\right.$$
$$\left.\times
\partial^\alpha_{x'}\partial^k_{x_n}\left[\partial^\delta_xf_2\partial^{j+1}_{\xi_n}
\partial^\delta_\xi
\sigma_{l+|\delta|}^{F_{\widehat{g}}}\right]|
_{x_n=0}\right\}d\xi_n\sigma(\xi')d^{n-1}x',\eqno(2.6)$$ \noindent
where the sum is taken over $
r-k-|\alpha|+l-j-1=-n,~~r,l\leq-1,~~|\alpha|\geq0.$ Then we have\\

 \noindent {\bf
Theorem 2.1}~~{\it  The form pair
$(\Omega_n(f_1,f_2)(\widehat{g}),\Omega_{n-1}(f_1,f_2)(\widehat{g}))$
only depends on $g$ and does not depend on the extension
$\widehat{g}$. It is a uniquely determined conformal invariant pair
on $[(M,g)]$ by (2.4), (2.5), and is symmetric in $f_1$ and
$f_2$.\\}

\noindent {\bf Proof.} By (2.5), (2.6), in order to prove that the
form pair
$(\Omega_n(f_1,f_2)(\widehat{g}),\Omega_{n-1}(f_1,f_2)(\widehat{g}))$
only depends on $g$ and does not depend on the extension
$\widehat{g}$, we only need to prove that
$D^{\alpha}_x(\sigma^{F_{\widehat{g}}}_{-j})|_M$ and
$D^{\alpha}_x(\sigma^{F_{\widehat{g}}}_{-j})|_{x_n=0}$ do not depend
on the extension $\widehat{g}.$ By Lemma A.3 in [U], this is
equivalent to proving that $D^{\alpha}_x(\widehat{g}_{i,j})|_M$ and
$D^{\alpha}_x(\widehat{g}_{i,j})|_{x_n=0}$ do not depend on the
extension $\widehat{g}$, where $[\widehat{g}_{i,j}]$ is the metric
matrix of $\widehat{g}$. The latter is trivial, so we prove the
first assertion. This fact says that
$(\Omega_n(f_1,f_2)(\widehat{g}),\Omega_{n-1}(f_1,f_2)(\widehat{g}))$
is a form pair with coefficients of derivatives of $g_{i,j}$, so we
can write $(\Omega_n(f_1,f_2)(g),\Omega_{n-1}(f_1,f_2)(g))$ instead
of
 $(\Omega_n(f_1,f_2)(\widehat{g}),\Omega_{n-1}(f_1,f_2)(\widehat{g})).$
\\
\indent By (2.5),(2.6), in order to prove
$(\Omega_n(f_1,f_2)(g),\Omega_{n-1}(f_1,f_2)(g))$ is a conformal
invariant of $[(M,g)]$, we only need prove
$\int_{|\xi|=1}p_{-n}(x,\xi)\sigma(\xi);~\int_{|\xi'|=1}p'_{-n+1}(x',\xi')\sigma(\xi')$
where $p_{-n}(x,\xi)~( p'_{-n+1}(x',\xi'))$ is a homogeneous
function of degree ${-n}~ (-n+1)$ about $\xi~ (\xi')$ and
$D^{\alpha}_x(\sigma^{F_{\widehat{g}}}_{-j})|_M$;~
$D^{\alpha}_x(\sigma^{F_{\widehat{g}}}_{-j})|_{x_n=0}$ do not depend
on the choice of the representative of $[(M,g)]$. As the discussions
in [AM],
$\int_{|\xi|=1}p_{-n}(x,\xi)\sigma(\xi);~\int_{|\xi'|=1}p'_{-n+1}(x',\xi')\sigma(\xi')$
do not depend on the choice of metric. For any representative $e^fg$
of $[(M,g)]$ where $f\in C^{\infty}(M)$, since
$(\Omega_n(f_1,f_2)(g),\Omega_{n-1}(f_1,f_2)(g))$ does not depend on
the extension $\widehat{g}$, we can choose the extension
$e^{\widehat{f}}\widehat{g}$ of $e^fg$ to compute
$(\Omega_n(f_1,f_2)(e^fg),$\\
\noindent $\Omega_{n-1}(f_1,f_2)(e^fg))$ where $\widehat{f}\in
C^{\infty}(\widehat{M})$ is an extension of $f$. By
$F_{\widehat{g}}=F_{e^{\widehat{f}}\widehat{g}}$, so symbols
$\sigma(F_{\widehat{g}})=\sigma(F_{e^{\widehat{f}}\widehat{g}})$.
Then by (2.5) and (2.6),
$$(\Omega_n(f_1,f_2)(g),\Omega_{n-1}(f_1,f_2)(g))
=(\Omega_n(f_1,f_2)(e^fg),\Omega_{n-1}(f_1,f_2)(e^fg)).$$ The other
properties of $(\Omega_n(f_1,f_2)(g),\Omega_{n-1}(f_1,f_2)(g))$
come from Theorem 3.1 and Proposition 3.3 in [Wa2]. \hfill$\Box$\\

\section{The Computation of
$(\Omega_4(f_1,f_2),\Omega_3(f_1,f_2))$}

 \quad In this section, we want to compute $(\Omega_4(f_1,f_2),\Omega_3(f_1,f_2))$ defined in Section 2 when $n=4$. We hope
 to compare the change of $(\Omega_4(f_1,f_2),\Omega_3(f_1,f_2))$ under the product metric and the nonproduct metric. So for simplicity,
  we firstly assume that $(\star)$ $f_1,~f_2$ {\it are independent of $x_n$ near the boundary.} For the general case, we will point out
  the way of computations in Section 4.\\
\indent $\Omega_4(f_1,f_2)$ is computed by Theorem 4.5 in [Wa2]. By
(2.6) and the assumption $(\star)$, then
$$\Omega_3(f_1,f_2)=\sum^{\infty}_{j,k=0}\sum
\sum^{-r}_{|\beta'|=1}\sum^{-l}_{|\delta'|=1}
\frac{(-i)^{j+k+1+|\alpha|+|\beta'|+|\delta'|}}{\alpha!\beta'!\delta'!(j+k+1)!}$$
$$\times
\int_{|\xi'|=1} \int^{+\infty}_{-\infty}{\rm
trace}_{\wedge^2T^*M}\left\{
\left[\partial^{\beta'}_{x'}f_1\times\partial^j_{x_n}\partial^{\alpha+\beta'}_{\xi'}\partial^k_{\xi_n}\pi^+_{\xi_n}
\sigma_{r+|\beta'|}^{F_{\widehat{g}}}\right]|_{x_n=0}\right.$$
$$\left.\times
\partial^\alpha_{x'}\left[\partial^{\delta'}_{x'}f_2\partial^k_{x_n}\partial^{j+1}_{\xi_n}
\partial^{\delta'}_{\xi'}
\sigma_{l+|\delta'|}^{F_{\widehat{g}}}\right]|
_{x_n=0}\right\}d\xi_n\sigma(\xi')d^{n-1}x',\eqno(3.1)$$ \noindent
where the sum is taken over $
-(r+l)+|\alpha|+k+j=3,~~r,l\leq-1,~~\alpha,\beta',\delta' \in {\bf
Z}^3_+.$ Since $\Omega_3(f_1,f_2)$ is a global form on $\partial M$,
so for any fixed point $x_0\in\partial M$, we can choose the normal
coordinates $V$ of $x_0$ in $\partial M$ (not in $M$) and compute
$\Omega_3(f_1,f_2)(x_0)$ in the coordinates
$x=(x',x_n)=(x_1,\cdots,x_{n-1},x_n)$ and domain
$\widetilde{V}=V\times [0,1)\subset M$ and the metric
$\frac{1}{h(x_n)}g^{\partial M}+dx_n^2.$ The dual metric of $g^M$ on
$\widetilde{V}$ is ${h(x_n)}g^{\partial M}+dx_n^2.$ Write
$g^M_{ij}=g^M(\frac{\partial}{\partial x_i},\frac{\partial}{\partial
x_j});~ g_M^{ij}=g^M(dx_i,dx_j)$, then
$$[g^M_{i,j}]= \left[\begin{array}{lcr}
  \frac{1}{h(x_n)}[g_{i,j}^{\partial M}]  & 0  \\
   0  &  1
\end{array}\right];~~~
[g_M^{i,j}]= \left[\begin{array}{lcr}
  h(x_n)[g^{i,j}_{\partial M}]  & 0  \\
   0  &  1
\end{array}\right],$$
and
$$\partial_{x_s}g_{ij}^{\partial M}(x_0)=0, 1\leq i,j\leq
n-1; ~~~g_{ij}^M(x_0)=\delta_{ij}.\eqno(3.2)$$
 \noindent We'll compute ${\rm
tr}_{\wedge^2(T^*M)}$ in the frame $\{dx_{i_1}\wedge dx_{i_2}|~1\leq
i_1<i_2\leq 4\},$ which is independent of the choice of frames.
 Let $\epsilon (\xi),~\iota (\xi)$ be the exterior and interior multiplications
 respectively where $\xi=\sum_{i=1}^n\xi_idx_i$ denotes a cotangent vector.
Recall Lemma 2.2 in [Wa1]
$$\partial_{x_j}(|\xi|_{g^M}^2)(x_0)=0,~{\rm if
}~j<n;~\partial_{x_n}(|\xi|_{g^M}^2)(x_0)=h'(0)|\xi'|_{g^{\partial
M}}^2.\eqno(3.3)$$ \noindent By (3.2) and $h(0)=1$, then under the
frame $\{dx_{i_1}\wedge dx_{i_2}|~1\leq i_1<i_2\leq 4\},$~\\
\noindent $\partial_{x_i}\varepsilon(dx_j)=0$ and
$$\partial_{x_l}\iota(dx_j)(x_0)=0,~{\rm if
}~l<n;~\partial_{x_n}\iota(dx_j)(x_0)=h'(0)\iota(dx_j)(x_0).\eqno(3.4)$$
\noindent So if $i<n$, then
$$\partial_{x_i}\varepsilon(\xi)(x_0)=\partial_{x_i}\iota(\xi)(x_0)=0;~
\partial_{x_n}\iota(\xi)(x_0)=h'(0)\iota(\xi')(x_0).\eqno(3.5)$$

\noindent {\bf Theorem 3.1}~~{\it Under the above conditions,}
$$\Omega_3(f_1,f_2)(x_0)=h'(0)\sum_{1\leq
i,j\leq3}a_{i,j}\partial_{ x_i}f_1\partial_{x_j}f_2dx_1\wedge
dx_2\wedge dx_3, \eqno(3.6)$$ \noindent {\it where $a_{i,j}$ is a
constant.}\\

\noindent {\bf Corollary 3.2}~~{\it Under the assumption $(\star)$,
if $h'(0)=0$ (for example $h=1-x_n^2$), then $\Omega_3(f_1,f_2)=0$.
Especially, if $g^{M}$ has the product metric near the boundary,
then $\Omega_3(f_1,f_2)=0$ and}
$$\widetilde{{\rm
Wres}}({\pi^+f_0}[\pi^+F_{\widehat{g}},\pi^+f_1][\pi^+F_{\widehat{g}},\pi^+f_2])=\int_{M}f_0\Omega
_{4}(f_1,f_2).\eqno(3.7)
$$

\indent Now we prove Theorem 3.1. Since the sum is taken over $
-(r+l)+|\alpha|+k+j=3,~~r,l\leq-1$, so $\Omega_3(f_1,f_2)$ is the
sum of the following five cases.\\

\noindent  {\bf case a)~I)}~$r=-1,~l=-1~k=j=0,~|\alpha|=1$\\

\noindent For convenience, we use $F$ instead of $F_{\widehat{g}}$
in the following. Let $\sigma_{L}^{F}$ denote the leading symbol of
$F$. By (3.1), we get
$${\rm case~a)~I)}=
\sum_{|\alpha|=1}\sum_{|\beta'|=1}\sum_{|\delta'|=1}
 \int_{|\xi'|=1}
\int^{+\infty}_{-\infty}{\rm trace}_{\wedge^2T^*M}\left\{
\partial^{\beta'}_{x'}f_1\partial^{\alpha+\beta'}_{\xi'}\pi^+_{\xi_n}
\sigma_{L}^{F}\right.$$
$$\left.\times
\left[\partial^{\alpha+\delta'}_{x'}f_2\partial_{\xi_n}
\partial^{\delta'}_{\xi'}
\sigma_{L}^{F}+\partial^{\delta'}_{x'}f_2\partial^\alpha_{x'}\partial_{\xi_n}
\partial^{\delta'}_{\xi'}
\sigma_{L}^{F}\right]\right\}(x_0)
d\xi_n\sigma(\xi')d^{n-1}x'.\eqno(3.8)$$ \noindent It is necessary
to compute
$${\rm trace}_{\wedge^2T^*M}[
\partial^{\alpha+\beta'}_{\xi'}\pi^+_{\xi_n}
\sigma_{L}^{F} \times
\partial_{\xi_n}
\partial^{\delta'}_{\xi'}
\sigma_{L}^{F}](x_0)$$ \noindent and  $${\rm trace}_{\wedge^2T^*M}[
\partial^{\alpha+\beta'}_{\xi'}\pi^+_{\xi_n}
\sigma_{L}^{F} \times
\partial_{\xi_n}
\partial^{\delta'}_{\xi'}
\partial^\alpha_{x'}\sigma_{L}^{F}](x_0).$$ \noindent Using the computations in
[Wa2,p.17], for $l,i,j<n$, then
$$\partial_{\xi_l}\partial_{\xi_i}\partial_{\eta_j}\left\{{\rm
trace}\left[\pi^+_{\xi_n}\sigma_L(F)(\xi',\xi_n)\times
\partial_{\xi_n}\sigma_L(F)(\eta',\xi_n)\right]\right\}|_{\xi'=\eta'}=\sum\xi_{i_1}\xi_{i_2}\cdots\xi_{i_{2k+1}}f(\xi_n),$$
\noindent where $f(\xi_n)$ is a smooth function about $\xi_n$ and
$1\leq i_1,\cdots,i_{2k+1}<n.$ Integration over $|\xi'|=1$ is zero.
By (3.3) and (3.5), then
$$\partial_{x_i}\sigma_{L}^{F}(x_0)=\partial_{x_i}\left[
\frac{\varepsilon(\xi)\iota(\xi)-\iota(\xi)\varepsilon(\xi)}{|\xi|^2}\right](x_0)=0,$$
\noindent so  case a)~I) yields zero.\\

\noindent  {\bf case a)~II)}~$r=-1,~l=-1~k=|\alpha|=0,~j=1$\\

By (3.1), we get
$${\rm case~a)~II)}=\frac{1}{2}
\sum_{|\beta'|=1}\sum_{|\delta'|=1}
 \int_{|\xi'|=1}
\int^{+\infty}_{-\infty}\partial^{\beta'}_{x'}f_1\partial^{\delta'}_{x'}f_2$$
$$\times{\rm
trace}_{\wedge^2T^*M}[
\partial^{\beta'}_{\xi'}\pi^+_{\xi_n}\partial_{x_n}
\sigma_{L}^{F}\times\partial^{\delta'}_{\xi'}
\partial^2_{\xi_n}
\sigma_{L}^{F}](x_0) d\xi_n\sigma(\xi')d^{n-1}x'.\eqno(3.9)$$

\noindent  {\bf case a)~III)}~$r=-1,~l=-1~j=|\alpha|=0,~k=1$\\

By (3.1), we get
$${\rm case~a)~III)}=\frac{1}{2}
\sum_{|\beta'|=1}\sum_{|\delta'|=1}
 \int_{|\xi'|=1}
\int^{+\infty}_{-\infty}\partial^{\beta'}_{x'}f_1\partial^{\delta'}_{x'}f_2$$
$$\times{\rm
trace}_{\wedge^2T^*M}[
\partial^{\beta'}_{\xi'}\partial_{\xi_n}\pi^+_{\xi_n}
\sigma_{L}^{F}\times\partial^{\delta'}_{\xi'}
\partial_{\xi_n}\partial_{x_n}
\sigma_{L}^{F}](x_0) d\xi_n\sigma(\xi')d^{n-1}x'.\eqno(3.10)$$
\noindent Write
$$p(\xi)=\varepsilon(\xi)\iota(\xi)-\iota(\xi)\varepsilon(\xi).$$
\noindent By (3.3),(3.4),(3.5), then
$$\partial_{x_n}p(\xi)(x_0)=h'(0)[\varepsilon(\xi)\iota(\xi')-\iota(\xi')\varepsilon(\xi)](x_0);$$
$$\partial_{x_n}\sigma_{L}^{F}(x_0)=\frac{h'(0)[\varepsilon(\xi)\iota(\xi')-\iota(\xi')\varepsilon(\xi)](x_0)}
{|\xi|^2}-\frac{h'(0)|\xi'|^2p(\xi)}{|\xi|^4}.
$$
\noindent So case a)~II+III) has the form in Theorem
3.1.\\

\noindent  {\bf case b)}~$r=-2,~l=-1,~k=j=|\alpha|=0$\\

By (3.1), we get
$${\rm case~b)}=
\sum_{|\beta'|=1}^2\sum_{|\delta'|=1}\frac{(-i)^{2+|\beta'|}}{\beta'!}
 \int_{|\xi'|=1}
\int^{+\infty}_{-\infty}\partial^{\beta'}_{x'}f_1\partial^{\delta'}_{x'}f_2$$
$$\times{\rm
trace}_{\wedge^2T^*M}[
\partial^{\beta'}_{\xi'}\pi^+_{\xi_n}
\sigma_{-2+|\beta'|}^{F}\times\partial^{\delta'}_{\xi'}
\partial_{\xi_n}
\sigma_{L}^{F}](x_0) d\xi_n\sigma(\xi')d^{n-1}x'.\eqno(3.11)$$
\noindent When $|\beta'|=2$, the term
$$\partial_{\xi_l}\partial_{\xi_i}\partial_{\eta_j}\left\{{\rm
trace}\left[\pi^+_{\xi_n}\sigma_L(F)(\xi',\xi_n)\times
\partial_{\xi_n}\sigma_L(F)(\eta',\xi_n)\right]\right\}|_{\xi'=\eta'}$$
will appear, as the disscusions in line 4 on p.6, it is zero after
the integration over $|\xi'|=1$.  So $|\beta'|=1$. In the following,
we prove that $\sigma_{-1}(F)(x_0)$ has the coefficient $h'(0)$.
Write $F=\frac{A}{\triangle},$ where $A=d\delta-\delta
d,~\triangle=d\delta+\delta d,$ then by the composition formula of
the symbol, we have
$$\sigma(F)=\sum_{|\alpha|\geq0}\sum_{0\leq i\leq 2}\sum_{j\geq
2}\frac{1}{\alpha!}\partial_{\xi}^{\alpha}(\sigma_i(A))D^{\alpha}_x(\sigma_{-j}(\triangle^{-1}));$$
$$\sigma_{-1}(F)=\sigma_1(A)\sigma_{-2}(\triangle^{-1})+\sigma_2(A)\sigma_{-3}(\triangle^{-1})+\sum_{|\alpha|=1}
 \partial_{\xi}^{\alpha}(\sigma_2(A))D^{\alpha}_x(\sigma_{-2}(\triangle^{-1})).\eqno(3.12)$$
\noindent By (3.3), then $$\sum_{|\alpha|=1}
 \partial_{\xi}^{\alpha}(\sigma_2(A))D^{\alpha}_x(\sigma_{-2}(\triangle^{-1}))(x_0)
 =\frac{ih'(0)|\xi'|^2\partial_{\xi_n}p(\xi)}{|\xi|^4}.\eqno(3.13)$$
\noindent Similar to (3.12), then
$$\sigma_1(d\delta)=\sigma_1(d)\sigma_0(\delta)+\sigma_0(d)\sigma_1(\delta)
-\sqrt{-1}\sum_i\partial_{\xi_i}(\sigma_1(d))\partial_{x_i}(\sigma_1(\delta));$$
$$\sigma_1(\delta
d)=\sigma_1(\delta)\sigma_0(d)+\sigma_0(\delta)\sigma_1(d)
-\sqrt{-1}\sum_i\partial_{\xi_i}(\sigma_1(\delta))\partial_{x_i}(\sigma_1(d)).\eqno(3.14)$$
\indent Let $\{e_1,\cdots,e_{n-1}\}$ be the orthonormal frame field
in $V$ about $g^{\partial M}$ which is parallel along geodesics and
$e_i(x_0)=\frac{\partial}{\partial x_i}(x_0)$, then
$\{\widetilde{e_1}=\sqrt{h(x_n)}e_1,\cdots,\widetilde{e_{n-1}}=\sqrt{h(x_n)}e_{n-1},
\widetilde{e_n}=dx_n\}$ is the orthonormal frame field in
$\widetilde V$ about $g^M$. By Lemma 2.3 and Section 3 in [Wa1], we
have
$$\sigma_1(d)=\sqrt{-1}\varepsilon(\xi),~\sigma_0(d)(x_0)=\frac{1}{4}h'(0)\sum^{n-1}_{i=1}\varepsilon(e_i^*)[\overline{c}(e_n)\overline{c}(e_i)
-c(e_n)c(e_i)];\eqno(3.15)$$
$$\sigma_1(\delta)=-\sqrt{-1}\iota(\xi),~\sigma_0(\delta)(x_0)=-\frac{1}{4}h'(0)\sum^{n-1}_{i=1}\iota(e_i^*)[\overline{c}(e_n)\overline{c}(e_i)
-c(e_n)c(e_i)],\eqno(3.16)$$ \noindent where
$$c(e_j)=\varepsilon(e_j^*)-\iota(e_j^*),~\overline{c}(e_j)=\varepsilon(e_j^*)+\iota(e_j^*).$$
\noindent By (3.5), then
$$\sigma_1(d\delta)(x_0)=\sqrt{-1}\varepsilon(\xi)\sigma_0(\delta)(x_0)-\sqrt{-1}\sigma_0(d)(x_0)\iota(\xi)
-\sqrt{-1}h'(0)\varepsilon(dx_n)\iota(\xi')(x_0);(3.17)$$
$$\sigma_1(\delta
d)(x_0)=-\sqrt{-1}\iota(\xi)\sigma_0(d)(x_0)+\sqrt{-1}\sigma_0(\delta)(x_0)\varepsilon(\xi).\eqno(3.18)$$
\noindent By Lemma A.1 in [U] and (3.3), then
$$\sigma_{-3}(\triangle^{-1})(x_0)=-\frac{1}{|\xi|^2}[\sigma_1(\triangle)\frac{1}{|\xi|^2}-
\sqrt{-1}\sum_i\partial_{\xi_i}(|\xi|^2)\partial_{x_i}(\frac{1}{|\xi|^2})](x_0)$$
$$~~~~~~~~~~~~~=-\frac{\sigma_1(\triangle)(x_0)}{|\xi|^4}-\frac{2\sqrt{-1}h'(0)|\xi'|^2\xi_n}{|\xi|^6}.
\eqno(3.19)$$ \noindent By (3.12), (3.13), (3.15)-(3.19) and the
definitions of $A,~\triangle$, we get
$\sigma_{-1}(F)(x_0)=h'(0)f(\xi).$ So case
b) has the form in Theorem 3.1.\\

\noindent {\bf  case c)}~$r=-1,~l=-2,~k=j=|\alpha|=0$\\

By (3.1), we get
$${\rm case~b)}=
\sum_{|\beta'|=1}\sum^2_{|\delta'|=1}\frac{(-i)^{2+|\delta'|}}{\delta'!}
 \int_{|\xi'|=1}
\int^{+\infty}_{-\infty}\partial^{\beta'}_{x'}f_1\partial^{\delta'}_{x'}f_2$$
$$\times{\rm
trace}_{\wedge^2T^*M}[
\partial^{\beta'}_{\xi'}\pi^+_{\xi_n}
\sigma_L^{F}\times\partial^{\delta'}_{\xi'}
\partial_{\xi_n}
\sigma_{-2+|\delta'|}^{F}](x_0)
d\xi_n\sigma(\xi')d^{n-1}x'.\eqno(3.20)$$ \noindent Similar to the
discussions in case b), case c) also has the form in Theorem
3.1, so we proved Theorem 3.1. \hfill$\Box$\\

\section{Some Remarks}

\quad In this section, we will point out the way of computations of
$a_{ij}$ in Theorem 3.1 and $\Omega_3(f_1,f_2)$ in the case of
$f_1,f_2$ depending on $x_n$ by some remarks.\\

 \noindent {\bf Remark 1} Since the computation of $\pi^+_{\xi_n}
\sigma_{-1}^{F}(x_0)$ is a little tedious, the computation of case
c) is more direct than the computation of case b). So we try to use
the computation of case c) and some simple computations instead of
the computation of case b). By the Leibniz rule, trace property and
"++" and "-~-" vanishing after the integration over $\xi_n$ (for
details, see [FGLS]), then
\begin{eqnarray*}
& &\int^{+\infty}_{-\infty}{\rm trace}_{\wedge^2T^*M}[
\partial_{\xi_n}\pi^+_{\xi_n} \sigma_{-1}^{F}(\xi',\xi_n)\times
\sigma_L^{F}(\eta',\xi_n)](x_0)d\xi_n\\
&=&\int^{+\infty}_{-\infty}{\rm trace}_{\wedge^2T^*M}[
\partial_{\xi_n} \sigma_{-1}^{F}(\xi',\xi_n)\times
\sigma_L^{F}(\eta',\xi_n)](x_0)d\xi_n\\
~~~~~~~~~~~~&&-\int^{+\infty}_{-\infty}{\rm trace}_{\wedge^2T^*M}[
\partial_{\xi_n}\pi^-_{\xi_n} \sigma_{-1}^{F}(\xi',\xi_n)\times
\sigma_L^{F}(\eta',\xi_n)](x_0)d\xi_n\\
&=&-\int^{+\infty}_{-\infty}{\rm trace}_{\wedge^2T^*M}[
 \sigma_{-1}^{F}(\xi',\xi_n)\times
\partial_{\xi_n}\sigma_L^{F}(\eta',\xi_n)](x_0)d\xi_n\\
~~~~~~~~~~~~~~&&- \int^{+\infty}_{-\infty}{\rm
trace}_{\wedge^2T^*M}[
\pi^+_{\xi_n}\sigma_L^{F}(\eta',\xi_n)\times\partial_{\xi_n}\pi^-_{\xi_n}
\sigma_{-1}^{F}(\xi',\xi_n)](x_0)d\xi_n\\
&=&-\int^{+\infty}_{-\infty}{\rm trace}_{\wedge^2T^*M}[
 \sigma_{-1}^{F}(\xi',\xi_n)\times
\partial_{\xi_n}\sigma_L^{F}(\eta',\xi_n)](x_0)d\xi_n\\
~~~~~~~~~~~~~&&- \int^{+\infty}_{-\infty}{\rm trace}_{\wedge^2T^*M}[
\pi^+_{\xi_n}\sigma_L^{F}(\eta',\xi_n)\times\partial_{\xi_n}
\sigma_{-1}^{F}(\xi',\xi_n)](x_0)d\xi_n.
\end{eqnarray*}
\noindent For computations of case a) II) and III), we have a
similar remark. But we may not get the sum of case b) and case c) is
zero through the above computations although we conjecture that it
should vanish and $\Omega_3(f_1,f_2)$
is also zero.\\

\noindent {\bf Remark 2} The computations of the trace of some
operators will appear in this case. We just compute an example and
the others are similar. In the following, we compute the equality:
$${\rm trace}_{\wedge^2T^*M}\{[\partial
_{x_n}p(\xi)]p(\eta)\}(x_0)=h'(0)[a_{n,m}\langle\xi',\eta'\rangle^2+b_{n,m}|\xi'|^2|\eta|^2](x_0)
+8h'(0)\xi_n\eta_n\langle\xi',\eta'\rangle,\eqno(4.1)$$ \noindent
where $C^m_n-a_{n,m}=b_{n,m}=C^{m-2}_{n-2}+C^m_{n-2}-2C^{m-1}_{n-2}$
and $C^m_n=\frac{n!}{m!(n-m)!}.$\\

\noindent {\bf Proof.}~By (3.5),
$$\partial_{x_n}
p(\xi)(x_0)=h'(0)[\varepsilon(\xi)\iota(\xi')-\iota(\xi')\varepsilon(\xi)](x_0)=h'(0)p(\xi',0)+\xi_nB,\eqno(4.2)$$
\noindent where
$B=h'(0)[\varepsilon(dx_n)\iota(\xi')-\iota(\xi')\varepsilon(dx_n)](x_0).$
By the well-known equality
$$\varepsilon_{m-1}(\xi)\iota_m(\eta)+\iota_{m+1}(\eta)\varepsilon_m(\xi)=\langle\xi,\eta\rangle
I_m,\eqno(4.3)$$ \noindent then
$$\varepsilon(dx_n)\iota(\xi')-\iota(\xi')\varepsilon(dx_n)=2\varepsilon(dx_n)\iota(\xi');~
p(\eta)=2\varepsilon(\eta)\iota(\eta)-\langle\eta,\eta\rangle
I_m.\eqno(4.4)$$ \noindent So by (4.2), (4.4) and Theorem 4.3 in
[U],
$${\rm trace}_{\wedge^2T^*M}\{[\partial
_{x_n}p(\xi)]p(\eta)\}(x_0)=h'(0)[a_{n,m}\langle\xi',\eta'\rangle^2+b_{n,m}|\xi'|^2|\eta|^2](x_0)$$
$$~~~+4h'(0)\xi_n{\rm
trace}_{\wedge^2T^*M}[\varepsilon(dx_n)\iota(\xi')\varepsilon(\eta)\iota(\eta)]-2|\eta|^2h'(0)\xi_n
{\rm trace}_{\wedge^2T^*M}[\varepsilon(dx_n)\iota(\xi')]\eqno(4.5)$$
\noindent By (4.3) and the trace property, we have
$${\rm
trace}_{\wedge^2T^*M}[\varepsilon(dx_n)\iota(\xi')]=0.\eqno(4.6)$$
\noindent As in [U,p.12-13], we write
$$a_m(\xi_1,\xi_2,\eta_1,\eta_2)={\rm
trace}_{\wedge^mT^*M}[\varepsilon_{m-1}(\xi_1)\iota_m(\xi_2)\varepsilon_{m-1}(\eta_1)\iota_m(\eta_2)].$$
\noindent then
$$a_{m+1}(\eta_1,\xi_2,\xi_1,\eta_2)=a_m(\xi_1,\xi_2,\eta_1,\eta_2)+\langle\xi_1,\xi_2\rangle
\langle\eta_1,\eta_2\rangle[2A_{n,m}-C_n^m],\eqno(4.7)$$ \noindent
where $A_{n,m}=C_n^m-C_n^{m-1}+\cdots+(-1)^mC_n^0.$ So
$a_1(\xi_1,\xi_2,\eta_1,\eta_2)=\langle\eta_2,\xi_1\rangle
\langle\xi_2,\eta_1\rangle$ and
$$a_2(\eta_1,\xi_2,\xi_1,\eta_2)=\langle\eta_2,\xi_1\rangle
\langle\xi_2,\eta_1\rangle+\langle\xi_1,\xi_2\rangle
\langle\eta_1,\eta_2\rangle[2A_{n,1}-C_n^1].\eqno(4.8)$$ \noindent
So by (4.8) and $n=4$
$${\rm
trace}_{\wedge^2T^*M}[\varepsilon(dx_n)\iota(\xi')\varepsilon(\eta)\iota(\eta)]=a_2(dx^n,\xi',\eta,\eta)=2\eta_n\langle\xi',\eta'\rangle.
\eqno(4.9)$$ \noindent By (4.5),(4.6) and (4.9), we prove the
equality (4.1).\hfill$\Box$\\

\noindent {\bf Remark 3} When $n=4$ and $f_1,f_2$ depend on $x_n$,
by (2.6) and considering the sum is taken over $
-(r+l)+|\alpha|+k+j=3,~~r,l\leq-1,~~1\leq
|\beta|=|\beta'|+\beta''\leq -r,~~1\leq
|\delta|=|\delta'|+\delta''\leq -l,$ similar to Section 3, we
compute $\Omega_{n-1}(f_1,f_2)(x_0)$ as the sum of $24$ cases about
$(r,l,k,j,\alpha,\beta',\beta'',\delta',\delta'')$. This can not add
to new technical difficulties except for a little tedious computations.\\

 \noindent{\bf
Acknowledgement:}~~The author is indebted to Professors Kefeng Liu
and Hongwei Xu for their help and hospitality. He thanks Professor
Weiping Zhang for his encouragement and support. He also thanks
the referee for his careful reading and helpful comments.\\

\noindent{\bf References}\\

\noindent [Ad] M. Adler, {\it On a trace functional for formal
pseudo-differential operators and the sympletic structure of the
Korteweg de Vries type equations}, Invent. Math. 50: 219-248, 1979.\\
\noindent [AM] P. M. Alberti and R. Matthes, {\it Connes' trace
formula and Dirac realization of Maxwell and Yang-Mills action},
Noncommutative geometry and the standard model of elementary
particle physics (Hesselberg,1999), 40-74,
Lecture Notes in Phys., 596, Springer, Berlin.\\
\noindent [C1] A. Connes, {\it The action functinal in
noncommutative
geometry}, Comm. Math. Phys., 117:673-683, 1998.\\
\noindent [C2] A. Connes, {\it Quantized calculus and applications},
XIth International Congress of Mathematical Physics (Paris,1994),
15-36, Internat Press, Cambridge, MA, 1995.\\
\noindent [FGLS] B. V. Fedosov, F. Golse, E. Leichtnam, and E.
Schrohe. {\it The noncommutative residue for manifolds with
boundary}, J. Funct.
Anal, 142:1-31,1996.\\
\noindent [Gr] G. Grubb, {\it Functional calculus for boundary value
problem}, Number 65 in Progress in Mathematics, Birkh\"{a}user,
Basel, 1986.\\
\noindent [Gu] V.W. Guillemin, {\it A new proof of Weyl's formula on
the asymptotic distribution of eigenvalues}, Adv. Math. 55 no.2,
131-160, 1985.\\
\noindent [K] D. Kastler, {\it The Dirac operator and gravitation},
Commun. Math. Phys, 166:633-643, 1995.\\
\noindent [KW] W. Kalau and M.Walze, {\it Gravity, non-commutative
geometry, and the Wodzicki residue}, J. Geom. Phys., 16:327-344, 1995.\\
\noindent [M] Yu. I. Manin, {\it Algebraic aspects of nonlinear
differential equations}, J. Sov. Math. 11: 1-22. 1979.\\
\noindent [S] E. Schrohe, {\it Noncommutative residue, Dixmier's
trace, and heat trace expansions on manifolds with boundary}, Contemp. Math. 242, 161-186, 1999.\\
 \noindent [U] W. J. Ugalde, {\it Differential forms and the
Wodzicki residue}, arXiv: Math, DG/0211361.\\
\noindent [Wa1] Y. Wang, {\it Gravity and the Wodzicki residue for
manifolds with boundary}, preprint, available online at {\it
www.cms.zju.edu.cn/frontindex.asp?version=english/priprint} \\
 \noindent [Wa2] Y. Wang, {\it Differential forms
and the Wodzicki residue for manifolds with boundary}, to appear J.
Geom. Phys., available online
at {\it www.sciencedirect.com}\\
\noindent [Wo] M. Wodzicki,  {\it Local invariants of spectral
asymmetry}, Invent.Math. 75 no.1 143-178, 1984.\\

\end{document}